\documentclass[12pt]{article}
\usepackage{amsmath}
\usepackage{amsthm}
\usepackage{amssymb}
\usepackage{graphicx}
\usepackage{booktabs}
\usepackage{threeparttable}
\title{Random potentials for pinning models with Laplacian interactions}
\author{Chien-Hao Huang \\
Academia Sinica}
\date{ }

\newtheorem{theorem}{\bf Theorem}[section]
\newtheorem{lemma}[theorem]{\bf Lemma}
\newtheorem{corollary}[theorem]{\bf Corollary}

\newtheorem{proposition}[theorem]{\bf Proposition}

\begin{document}
\maketitle
\fontsize{12}{18pt}\selectfont

\begin{abstract}
We consider a statistical mechanics model for biopolymers. Sophisticated polymer chains, such as DNA, have stiffness when they stretch chains. The Laplacian interaction is used to describe the stiffness. Also, the surface between two media has an attraction force, and the force will pull the chain back to the surface. In this talk, we deal with the random potentials when the monomers interact with the random media. Although these models are different from the pinning models studied before, the result about the gap between the annealed critical point and the quenched critical point stays the same.
\end{abstract}

\section{Introduction} 

\subsection{The interpretation of the model}

\noindent
The (1+1)-directed walk model has the walk $(i,\varphi_i)$ drifting in the space $\{\mathbb{N}\cup 0\}\times \mathbb{R}$. This model is one special kind of ``self-avoiding'' walk models. ``Self-avoiding'' describes the phenomenon that the two particles of the polymer chain can not occupy the same site. Since the first coordinate of $(i,\varphi_i)$ is strictly increasing, the self-avoidance is satisfied. The Laplacian interactions describe the physical phenomenon that some bio-polymers such as DNA dislike to bend too much. In this paper, we discuss the ``pinning'' model which explains that the polymer chain favours certain point or interface between two media. Every time the path of the polymer reaches the interface, a non-negative reward will be given. 

The general pinning model was discussed in monographs \cite{Gia07,Gia11} which point out that the polynomial decay of the renewal distribution of the walk plays an important role. Let the random walk $(S_n)_{n\geq 0}$, starting at $0$, describe the path of the polymer under the law $P$ with independent and identically distributed (i.i.d.) increments $(S_n-S_{n-1})_{ n\geq 1 }$. For every $h\in\mathbb{R}$, we define the partition function $Z_n(h)$ and the free energy $f(h)$ by
\begin{equation}
Z_n(h):= E(e^{h \sum_{i=1}^n \mathbf{1}_{S_i=0}}\cdot\mathbf{1}_{S_n=0})
\end{equation}
\begin{equation}\label{energy}
f(h):= \lim_{n\rightarrow\infty} \frac{1}{n}\log Z_n(h)
\end{equation}
Define the renewal sequence $(\tau_n)_{n\geq 0}$ by
\begin{equation}
\tau_0 :=0, \;\;\; \tau_{n+1}= \inf\{ k> \tau_n \; : \; S_k =0\}
\end{equation} 
Suppose that there exist a number $\alpha\in(0,\infty)$ such that 
\begin{equation}
K(n):= P(\tau_1=n)\sim \frac{c_K}{n^{1+\alpha}},
\end{equation}
The polymer path receive a reward $e^{h}$ when it reaches the origin. \cite{Gia11} has the following theorem.

\begin{proposition}
\cite{Gia11} For every $h\in\mathbb{R}$, $f(h)$ exists and is continuous. And there exists a number $h_c$ such that $f(h)=0$ for $h\leq h_c$ and $f(h)>0$ for $h> h_c$. Moreover
\begin{equation}
f(h) \stackrel{h\searrow h_c}{\sim} C(K(\cdot)) \left\{ \begin{array}{ll}
(h-h_c)^{1/\alpha},            & \alpha \in (0,1) \\
-(h-h_c)/\log (h-h_c),     & \alpha =1, \\
(h-h_c),              & \alpha >1.
\end{array}\right.
\end{equation}
where
\begin{equation}
 C(K(\cdot)) = \left\{ \begin{array}{ll}
C(K(\cdot),\alpha),            & \alpha \in (0,1) \\
1/c_K,     & \alpha =1, \\
\sum_n K(n)/ \sum_n nK(n),              & \alpha >1.
\end{array}\right.
\end{equation}
\end{proposition}

In the physical literature, if the free energy $f(h)$ has a continuous $(k-1)$th-derivative at a point $h$, but its $k$th-derivative is not continuous at the same point $h$, we say that the model {\it has $k$th order phase transition at the critical point $h$}. For example, Proposition 1.1 indicates that there is a first order phase transition at the critical point $h_c$ when $\alpha >1$, and there is a second order phase transition at the critical point $h_c$ when $\alpha =1$.

In this paper, we discuss the case that the polymer is influenced by the random environment. That is,  let the random field $\varphi:\{0,1,...,N\}\rightarrow \mathbb{R}^d$ represent the position of the polymer path and $V(\cdot)$ be the Gaussian potential $\frac{|x|^2}{2}$. The law of the field is given by $\exp\{-\sum_iV(\nabla\varphi_i)\}$ where $\nabla$ is the discrete gradient, and by  $\exp\{-\sum_iV(\Delta\varphi_i)\}$ where $\Delta$ is the discrete Laplacian. For every potential $V(\cdot)$, a random charge is added as a factor: $\exp(\beta\omega_i)V(\cdot)$ with $\omega_i$ satisfying the standard normal distribution. The partition function now depends on the environment $\omega=\{\omega_i\}$, denoted by $Z_{n,\omega}$. The {\it quenched free energy} is defined in the way as \eqref{energy}. If we average the environment first, that is, $\mathbb{E}Z_{n,\omega}$, then 
\begin{equation}
f^a:=\lim_{n\rightarrow\infty} \frac{1}{n}\log \mathbb{E}Z_{n,\omega}
\end{equation} 
is called the {\it annealed free energy}. Sometimes, $\lim_{n\rightarrow\infty} \frac{1}{n}\log Z_{n,\omega} = \lim_{n\rightarrow\infty} \frac{1}{n}\mathbb{E}\log Z_{n,\omega}$ almost surely. In this case, the annealed free energy is less than the quenched free energy by Jensen's inequality.
If the quenched and annealed critical points for the quenched and annealed free energy, respectively, both exist,  a question which physicists are interested in is ``are they equal to each other?'' In the ``weak disordered regime'', the gap between the annealed and quenched critical points is positive only when the disorder is large enough. In contrast, the ``strong disordered regime'' means that the gap between the annealed and quenched critical points is positive even when the disorder is small. We will give the answer to this question in the rest of the paper.

\subsection{The (1+1)-dimensional pinning model with $\Delta$-interaction}

We consider a model for biopolymers with the Gaussian potential $V(x):= \frac{|x|^2}{2}$ in the random environment. Let $\varphi:\{0,1,...,N\}\rightarrow \mathbb{R}$ be the position of the polymer path. The polymer measure is given by $\exp\{-\sum_i V(\Delta\varphi_i)\}$ where $\Delta$ is the discrete Laplacian. For every Gaussian potential $V(\cdot)$, a random charge is added as a factor: $\exp(\beta\omega_i)V(\cdot)$ with $\omega_i$ satisfying the standard normal distribution. The interaction with the origin in the random path space is also considered. Each time the path touches the origin, a reward $\epsilon\geq 0$ is given.

The Hamiltonian $H_N(\varphi):=H_{0,N}(\varphi)$ is defined as 
\begin{equation}
H_{M,N}(\varphi):= \sum_{n=M}^{N-1}\limits V(\Delta\varphi_{n}),
\;\; \Delta\varphi_{n}:= (\varphi_{n+1} -\varphi_{n})+(\varphi_{n-1} -\varphi_{n}),
\end{equation}
with boundary conditions $\varphi(M-1)=\varphi(M)=\varphi(N-1)=\varphi(N)=0$, where $V(x)$ is called the potential with $\int_{\mathbb{R}}\exp(-V(x)) \;dx=1$.
and the random Hamiltonian is defined by
\begin{equation}
H_{M,N}(\varphi):= \sum_{n=M}^{N-1}\limits e^{\beta\omega_n} \left( \frac{|\Delta\varphi_{n}|^2}{2} -\log(\sqrt{2\pi}) \right)
\end{equation}
and the polymer measure is given by
\begin{equation}
P^{\beta,\epsilon}_{N,\omega}(d\varphi_{1},\dots,d\varphi_{N-2}):=  \frac{e^{-H_{N,\omega}(\varphi) }}{Z^{\beta,\epsilon}_{N,\omega}} 
\prod_{i=1}^{N-2}\limits  ( \epsilon \delta_{0} (d\varphi_{i} )  +d\varphi_{i})
\end{equation}
The partition function $Z_{N,\omega}^{\beta,\epsilon}$ is defined as the normalizing constant.

The polymer chain obtains a reward $\epsilon$ for touching the origin in $\mathbb{R}$. If the displacements of the two consecutive segments of the path, $\varphi_{n+1} -\varphi_{n}$ and $\varphi_{n} -\varphi_{n-1}$, have different sign (path goes up and then down or first goes down then up), the corresponding Hamiltonian (probability) is bigger (smaller) than the Hamiltonian (probability) when they have the same sign. This characterizes the stiffness of the path since the path measure penalizes the paths which bend with big angles.

The non-random case (i.e. $\beta=0$) was discussed in $\cite{CD08,CD09}$ for the general potential $V(x)$.
Let $f(\epsilon)$ denote the free energy for the non-random case,
\begin{equation}
f(\epsilon):=\lim_{N\rightarrow \infty} f_N(\epsilon); 
\;\; f_N(\epsilon):= \frac{1}{N}\log Z_{N}^{\epsilon}.
\end{equation}

In $\cite{CD08}$, it was proved that there exists a positive number $\epsilon_c$ such that $f(\epsilon)=0$ for $\epsilon\leq \epsilon_c$ and $f(\epsilon)>0$ for $\epsilon > \epsilon_c$. Moreover, the phase transition for the pinning model is exactly of second order. We use the renewal equation to prove

\begin{proposition} There exists a constant $c_1$ such that
\begin{equation}
f(\epsilon_c e^{\delta}) \stackrel{\delta\searrow 0}{\sim}  \frac{c_1 \delta}{-\log\delta}.
\end{equation}
\end{proposition}

\noindent
From this proposition, it's easy to get

\begin{corollary} The phase transition is exactly of second order.
\end{corollary}

\noindent
The proof Proposition 1.2 is given in Section 2.2. One can see that it is the case the rate of polynomial decay of the renewal distribution has the exponent $\alpha=2$.

It's easy to see that, since we perturb every potential $V(\cdot)$, $f$ and $f^a$ are different when the randomness occur, namely, $\beta>0$. Thus, it is not interesting to consider the difference between the annealed and quenched critical points. So we introduce the ``adjusted'' quenched free energy
\begin{equation}F(\beta,\epsilon):=\lim_{N\rightarrow \infty} F_N(\beta,\epsilon); 
\;\;\; F_N(\beta,\epsilon):= \frac{1}{N}\log \frac{ Z_{N,\omega}^{\beta,\epsilon} }{Z_{N,\omega}^{\beta,0}},
\end{equation}
and the ``adjusted'' annealed free energy

\begin{equation}
F^a(\beta,\epsilon):=\lim_{N\rightarrow \infty} F^a_N(\beta,\epsilon); 
\;\;\; F^a_N(\beta,\epsilon):= \frac{1}{N}\log \mathbb{E} \frac{ Z_{N,\omega}^{\beta,\epsilon} }{Z_{N,\omega}^{\beta,0}}.
\end{equation}

Also, we introduce the ``adjusted'' free energy
\begin{equation}
F(\beta,\epsilon):=\lim_{N\rightarrow \infty} F_N(\beta,\epsilon); 
\;\; F_N(\beta,\epsilon):= \frac{1}{N}\log \frac{ Z_{N,\omega}^{\beta,\epsilon} }{Z_{N,\omega}^{\beta,0}}.
\end{equation}
The existence of the free energy will be proved in Section 2.3.

The following proposition states our main result of this paper.

\begin{proposition} 
Consider the ``adjusted'' free energy. There is a positive number $\beta_2$ such that for all $ 0< \beta < \beta_2$, the anneal critical point is strictly less than the quenched critical point.
\end{proposition}

The structure of this paper is the following: in Section 2, we first identify the asymptote for the Laplacian model in the non-random case, then we discuss the free energy and critical points. Some computations of matrix determinants are left in Section 3.

\section{The Laplacian model}

\noindent
In this section, we first consider the non-random case, which is discussed in \cite{CD08}.

\subsection{Known results for the non-random case}

\noindent
We define a contact process $(\tau_i)_{i\in\mathbb{N}\cup 0}$ by
\begin{equation}
\tau_0 :=0, \;\;\;   \tau_{i+1} :=  \inf\{  k>\tau_i \; :\; \varphi_k= 0\}
\end{equation}
\noindent
and the process $(J_i)_{i\in\mathbb{N}\cup 0}$, which is the height of the polymer path right before the path hits 0, namely, 
\begin{equation}
J_i:= \varphi_{\tau_i -1}
\end{equation}
and by definition $J_0 :=\varphi_{-1}=0$.
The quantity $L_s$ is the number of path contacts between $0$ and time $s$. 
\begin{equation}
L_s:= \# \{ 0 < n \leq s : \varphi_n =0\}
\end{equation}
The joint distribution of the process $\{ L_{N}, (\tau_i)_{i\leq L_{N}},  (J_i)_{i\leq L_{N}} \}$

$$P_{\epsilon, N} ( L_{N} =k , \tau_i=t_i, J_i\in dy_i , i=1,...,k-1)$$
\begin{equation}\label{3.0a}
=\frac{\epsilon^{k-1}}{Z_{\epsilon,\hat{N}}} F_{0,dy_1}(t_1)\cdot F_{y_1,dy_2}(t_2-t_1)
\cdots F_{y_{k-1},dy_k}(N-t_{k-1})\cdot F_{y_k,{0}}(1)
\end{equation}
where $F_{x,dy}(n):= f_{x,y}(n)\mu(dy) , \; \mu(dy):= \delta_0(dy)+dy$ and

\begin{equation}
f_{x,y}(n):= \left\{\begin{array}{ll}
\exp(-\beta x^2/2) \mathbf{1}_{y=0} &, n=1\\
\exp(-\beta H_{[-1,2](x,0,y,0)}) \mathbf{1}_{y\neq 0} &, n=2\\
\exp(-\beta H_{[-1,n](x,0,\varphi_1,\dots, \varphi_{n-2},y,0)}) \mathbf{1}_{y\neq 0} \;\; d\varphi_1 \cdots d\varphi_{n-2} &, n\geq 3
\end{array}
\right.
\end{equation}

We define
\begin{equation}
K_{x,dy}^{\epsilon}(n):= \epsilon F_{x,dy}(n)e^{f(\epsilon)n} \frac{v_{\epsilon}(y)}{v_{\epsilon}(x)},
\end{equation}
where $f(\epsilon)$ and $v_{\epsilon}(x)$ are given in \cite{CD08}.

We recall the following theorem from \cite{CD08}.
\begin{proposition} \cite{CD08} 
For each $\epsilon > 0$, there exist a real number $f(\epsilon)\in[0,\infty)$ and a positive real function $v_{\epsilon}(x)$ such that
\begin{equation}
\int_{y\in\mathbb{R}}\sum_{n\in \mathbb{N}} K^{\epsilon}_{x,dy}(n)=1, \;\;\forall x\in\mathbb{R}.
\end{equation}
\end{proposition}

Based on this proposition, we can define a measure $\mathcal{P}_{\epsilon}$ such that under $\mathcal{P}_{\epsilon}$ the process $\{\tau_i, J_i\}_{i\in \mathbb{N}\cup 0}$ is a Markov chain with the transition kernel 
\begin{equation}
\mathcal{P}_{\epsilon}\left(   (\tau_{i+1}, J_{i+1})\in (\{ n\},dy)            |(\tau_i,J_i)=(m,x)     \right)
=K^{\epsilon}_{x,dy}(n-m).
\end{equation}

The following proposition characterizes the relationship between $P_{\epsilon,N}$ and $\mathcal{P}_{\epsilon}$.

\begin{proposition}\cite{CD08}
Let $\mathcal{A}_N := \{ \exists i\geq 0 : \tau_i=N-1, \tau_{i+1}=N\}$. Then $\forall\; N$, $\epsilon>0$, $k\leq N$, $(t_i)_i$ and $(y_i)_i$
\begin{equation}
P_{\epsilon,N}(\ell_N=k , \tau_i=t_i , J_i\in dy_i ,i\leq k)=\mathcal{P}_{\epsilon}( 
\ell_N=k , \tau_i=t_i , J_i\in dy_i ,i\leq k | \mathcal{A}_N)
\end{equation}
{\it and}
\begin{equation}\label{imp1}
Z_{\epsilon,N}= \frac{e^{f(\epsilon)N}}{\epsilon^2} \mathcal{P}_{\epsilon}(\mathcal{A}_N).
\end{equation}
\end{proposition}

We define a process to describe the double-return, that is, $(\chi_i)_{i\in\mathbb{N}\cup 0}$ is defined by
\begin{equation}
\chi_0:=0, \;\;\; \chi_{i+1}:= \inf\{ k> \chi_i \; : \; \varphi_{k-1}=\varphi_k=0\}.
\end{equation}

\begin{proposition} \cite{CD08}
The process $(\chi_i)_i$ under $\mathcal{P}_{\epsilon}$ is a classical renewal process. Furthermore, there exists $\alpha>0$ such that for each $\epsilon\in[\epsilon_c,\epsilon_c +\alpha]$ we have
\begin{equation}\label{imp2}
\mathcal{P}_{\epsilon}(\chi_1=n)\sim \frac{C_{\epsilon}}{n^2}\exp(-f(\epsilon)n)
\end{equation}
as $n\rightarrow\infty$, where $C_{\epsilon}\in (0,\infty)$ is a continuous function of $\epsilon$.
\end{proposition}

\subsection{Proof of Proposition 1.2}

\noindent 
Denote \begin{equation}
\check{Z}_n:=Z_n(\{\mbox{no double returns} \}), \; n\geq 3,
\end{equation} which means that only the paths of ``no double returns'' are considered. Set 
$\check{Z}_{0,1}^{\epsilon}=Z_{0,1}^{\epsilon}=\frac{1}{\sqrt{2\pi}}, \; \check{Z}_{0,2}^{\epsilon}=0, Z_{0,2}^{\epsilon}=\frac{1}{2\pi}$. For $n\geq 3$, the renewal equation
\begin{equation}
Z_{0,n}^{\epsilon}=\check{Z}_{0,n}^{\epsilon}+ 
\check{Z}_{0,1}^{\epsilon}\epsilon Z_{1,n}^{\epsilon}+
\sum_{\chi =3}^{n-2} \check{Z}_{0,\chi}^{\epsilon}\epsilon^2 Z_{\chi,n}^{\epsilon} +  \check{Z}_{0,n-1}^{\epsilon}\epsilon Z_{n-1,n}^{\epsilon}.
\end{equation}
Define 
\begin{equation}
u_1=\frac{Z_{0,1}^{\epsilon}}{\epsilon} x; \;\;u_n=Z_{0,n}^{\epsilon} x^n,\; n=2,3,... 
\end{equation}
\begin{equation}
a_1=\epsilon\check{Z}_{0,1}^{\epsilon} x; \;\;a_n=\epsilon^2 \check{Z}_{0,n}^{\epsilon} x^n,\; n=2,3,... \end{equation}
\begin{equation}b_n=\check{Z}_{0,n}^{\epsilon} x^n,\; n=1,2,3,... \end{equation}
Thus,
\begin{equation}u_n=b_n+\sum_{i=1}^{n-1}a_iu_{n-i}.\end{equation}

Suppose that $x^{\epsilon}$ is the solution of 
\begin{equation}\label{3.1a}
\sum_{n\geq 1} a_n = \epsilon\check{Z}_{0,1}^{\epsilon} x +\sum_{n\geq 3}  \epsilon^2 \check{Z}_{0,n}^{\epsilon} x^n = 1,
\end{equation} 
then by \cite{Fel66} section XIII.4,
\begin{equation}
\lim_{n\rightarrow \infty} u_n=\frac{\sum_{n\geq 1} b_n}{\sum_{n\geq 1} na_n}.
\end{equation}
We have
\begin{equation}\label{3.1}
f(\epsilon)=-\ln x^{\epsilon}.
\end{equation} 
Note that $x^{\epsilon_c}=1$. Since 
\begin{equation}
\frac{\check{Z}^{\epsilon}_{0,n}}{Z_{0,n}^{\epsilon}}= P_{\epsilon,n}(\chi_1=n)= \mathcal{P}_{\epsilon}(\chi_1=n|\mathcal{A}_n)=\frac{\mathcal{P}_{\epsilon}(\chi_1=n)}{\mathcal{P}_{\epsilon}(\mathcal{A}_n)}
\end{equation} 
From \eqref{imp1}, we have
\begin{equation}
\check{Z}^{\epsilon}_{0,n}= \mathcal{P}_{\epsilon}(\chi_1=n) \cdot \frac{e^{f(\epsilon)n}}{\epsilon^2}
\end{equation} 

Now, we choose $\epsilon$ as $\epsilon_c+\delta$. Thanks to \eqref{imp2},
\begin{equation}
\check{Z}_{0,n}^{\epsilon}   \sim \frac{C_{\epsilon}}{\epsilon^2 n^2} \; \mbox{as}\; n\rightarrow\infty .
\end{equation}

Apply \eqref{3.1a} to $\epsilon$ and $\epsilon_c$,
$$(1-x)\check{Z}_{0,1}+\sum_{n=3}^{\infty}\left( \epsilon_c\check{Z}_{0,n}^{\epsilon_c}-\epsilon\check{Z}_{0,n}^{\epsilon}x^n\right)
=\frac{1}{\epsilon_c}-\frac{1}{\epsilon}.$$
So that
$$(1-x) \left[\check{Z}_{0,1}+\sum_{n\geq 3} \epsilon\check{Z}_{0,n}^{\epsilon}(1+\cdots+x^{n-1})\right]
+\sum_{n=3}^{\infty}\left( \epsilon_c\check{Z}_{0,n}^{\epsilon_c}-\epsilon\check{Z}_{0,n}^{\epsilon}\right)
=\frac{\epsilon-\epsilon_c}{\epsilon_c\epsilon}.$$
Then we rewrite the above equality to be
\begin{eqnarray}\label{3.1b}
(1-x)\left[ \check{Z}_{0,1}-\epsilon\check{Z}_{0,1} +\epsilon \sum_{j=0}^{\infty}\left(\sum_{n=j+1}^{\infty}\check{Z}_{0,n}^{\epsilon}\right)x^j\right]\\\label{3.1c}
=(\epsilon-\epsilon_c)\left[ \frac{1}{\epsilon_c\epsilon} +\left(\frac{\sum_{n=3}^{\infty} \epsilon\check{Z}_{0,n}^{\epsilon}-\sum_{n=3}^{\infty}\epsilon_c\check{Z}_{0,n}^{\epsilon_c}}{\epsilon-\epsilon_c}\right)\right].
\end{eqnarray}

In \eqref{3.1b}, we know $$\sum_{j=0}^{\infty}\left(\sum_{n=j+1}^{\infty}\check{Z}_{0,n}^{\epsilon}\right)x^j\sim -\frac{C_{\epsilon_c}}{\epsilon_c^2}\log (1-x) 
\; \mbox{as} \; x\nearrow 1,$$ and in \eqref{3.1c}, 
$$\frac{\sum_{n=3}^{\infty} \epsilon\check{Z}_{0,n}^{\epsilon}-\sum_{n=3}^{\infty}\epsilon_c\check{Z}_{0,n}^{\epsilon_c}}{\epsilon-\epsilon_c}\rightarrow  \frac{d}{d\epsilon}\left.\left(\sum_{n=3}^{\infty} \epsilon\check{Z}_{0,n}^{\epsilon}\right)\right|_{\epsilon=\epsilon_c +}:= c_0(\epsilon_c)
\; \mbox{as} \; \epsilon\searrow \epsilon_c,$$
since $\sum_{n=3}^{N} \epsilon\check{Z}_{0,n}^{\epsilon}$ is a convex polynomial in $\epsilon$, and it converges pointwise, so $\sum_{n=3}^{\infty} \epsilon\check{Z}_{0,n}^{\epsilon}$ is convex in $\epsilon$, thus, the right-hand derivative exists. Therefore, \eqref{3.1b} is asymptotic to
\begin{equation}\label{3.1d}
 -(1-x)\frac{C_{\epsilon_c}}{\epsilon_c}\log (1-x) 
\end{equation}
as $x\nearrow 1$ and \eqref{3.1c} is asymptotic to
\begin{equation}\label{3.1e}
(\epsilon-\epsilon_c)\left[ \frac{1}{\epsilon_c^2}+ c_0(\epsilon_c)\right]
\end{equation}
as $\epsilon\searrow \epsilon_c$.
Combine \eqref{3.1},  \eqref{3.1d}, and \eqref{3.1e},
\begin{equation}f(\epsilon_c+\delta) \stackrel{\delta\searrow 0}{\sim} \frac{c_1}{\epsilon_c}\cdot \frac{\delta}{-\log\delta}, \end{equation}
where
\begin{equation}c_1 =\frac{1+\epsilon_c^2 c_0(\epsilon_c)}{C_{\epsilon_c}}.\end{equation}
\noindent
Or equivalently, \begin{equation}f(\epsilon_c e^{\delta}) \stackrel{\delta\searrow 0}{\sim}  \frac{c_1 \delta}{-\log\delta}.\end{equation}
The proof of Proposition 1.2 is complete.

\subsection{Free energy for the random case}

\subsubsection{The annealed free energy} 

\noindent
Let $e^{-V_{\beta}(x)} =\mathbb{E}\left(\frac{\exp(\frac{\beta }{2}\omega_0)}{\sqrt{2\pi}}e^{-(\exp(\beta\omega_0))x^2/2}\right)$. Thus, $F^{a}(\beta,\epsilon)=F^{V_{\beta}}(\epsilon)$ and $\epsilon_c^{a}(\beta)=\epsilon_c^{V_{\beta}}>0$. We claim that $F(\epsilon M(-\beta)^{-\frac{1}{2}})\leq F^{a}(\beta,\epsilon)\leq F(\epsilon M(\frac{\beta}{2})))$. Therefore, 
$$\frac{1}{M( \frac{\beta}{2})}\leq \frac{\epsilon_c^a(\beta)}{ \epsilon_c(0)} \leq \sqrt{M(-\beta)}$$

The annealed partition function has upper bound
\[\begin{array}{rcl}
\mathbb{E}\mathcal{Z}^{\beta,\epsilon}_{N,\omega}
&\leq &  Z_{N}^{0,\epsilon M (\frac{\beta}{2})} \cdot (M(\frac{\beta}{2}))^{2}
\end{array}\]

and lower bound
\[\begin{array}{rcl}
\mathbb{E}  \mathcal{Z}^{\beta,\epsilon}_{N,\omega}
& =  & \mathbb{E} \left(\sum_{l=0}^{N-1} \epsilon^l (\sqrt{2\pi})^{-l}\left[\sum_{|p|=l+2}\limits c^{\varphi}_{ p} \;\prod_{i=0}^N \exp(-p_i \beta\omega_i)\right]^{-1/2}   \right)  \\
&\geq &  \sum_{l=0}^{N-1} \epsilon^l (\sqrt{2\pi})^{-l}\left[\sum_{|p|=l+2}\limits c^{\varphi}_{ p} \;\prod_{i=0}^N  \mathbb{E} \exp(-p_i \beta\omega_i)\right]^{-1/2}    \\

&=    &  \sum_{l=0}^{N-1} \epsilon^l (\sqrt{2\pi})^{-l}\left[\sum_{|p|=l+2}\limits c^{\varphi}_{ p} \; M(-\beta)^{l+2} \right]^{-1/2}    \\
&=    & \sum_{l=0}^{N-1} M(-\beta)^{-1} \left(\epsilon M(-\beta)^{-\frac{1}{2}} \right)^l (\sqrt{2\pi})^{-l}\left[\sum_{|p|=l+2}\limits c^{\varphi}_{ p}\right]^{-1/2}\\
&=    & M(-\beta)^{-1} \left(\sum_{l=0}^{N-1} (\epsilon M(-\beta)^{-\frac{1}{2}})^l (\sqrt{2\pi})^{-l}[\mbox{det}(L^{\varphi})]^{-1/2}\right)\\
&=    &  \frac{M(-\beta)^{-1}}{(2\pi)^{-1/2}}\;  Z_N^{0,\epsilon M(-\beta)^{-\frac{1}{2}}}
\end{array}\]
This proves the claim.\\

\subsubsection{$\epsilon=0$} 

\noindent
When $\epsilon=0$, $Z^{\beta,0}_{N+1,\omega}=\int_{\mathbb{R}^{N-1}} \frac{1}{\sqrt{2\pi}^N} \mbox{exp}(-\frac{1}{2} \langle \varphi,L^{\omega}\varphi\rangle) \prod_{i=1}^{N-1}\limits d\varphi_{i}$, where $L^{\omega}$ is a symmetric $(N-1)\times(N-1)$ matrix. The upper triangle part is defined as following:
\[l^{\omega}_{ij}= \left \{ \begin{array}{ll}
\exp(\beta\omega_{i-1})+4\exp(\beta\omega_{i})+\exp(\beta\omega_{i+1}), & i=j\\
-2\exp(\beta\omega_{i})-2\exp(\beta\omega_{i+1}),                      & i=j-1 \\
\exp(\beta\omega_{i+1}),                                         & i\leq j-2\\
0,                    & i\leq j-3.
\end{array}\right.\]

For the determinant of $L^{\omega}$, we have the following lemma.

\begin{lemma}
\begin{equation}\mbox{det}(L^{\omega}) = \prod_{i=0}^{N} \exp(\beta \omega_i)\cdot \left[   \sum_{k=1}^N  \sum_{i=0}^{N-k} k^2\exp(-\beta \omega_i)\exp(-\beta \omega_{i+k})\right].
\end{equation}
\end{lemma}
The proof is left in Section 3. Note that when $\beta=0$, det($L^{\omega}$)=$\frac{1}{12}N(N+1)^2(N+2)$.
\noindent
Thus, 
\begin{equation}
Z_{N+1,\omega}^{\beta,0} =  \frac{1}{\sqrt{2\pi}} \cdot \prod_{i=0}^{N} \exp(-\frac{\beta}{2} \omega_i)\cdot \left[   \sum_{k=1}^N  \sum_{i=0}^{N-k} k^2\exp(-\beta \omega_i)\exp(-\beta \omega_{i+k})\right]^{-1/2} 
\end{equation}
\noindent
Let $T_N= \sum_{0\leq i < j \leq N}\exp(-\beta \omega_i)\exp(-\beta \omega_{j})$. The term in the bracket is bounded by $T_N$ and $N^2T_N$. Since $\lim_{N\rightarrow \infty} (N+1)^{-2}T_N =\frac{1}{2}M(-\beta)^2$ a.s., we have \begin{equation}
\lim_{N\rightarrow \infty} \frac{1}{N}\log Z_{N,\omega}^{\beta,0} =0 \; \mbox{almost surely}.
\end{equation}

\subsubsection{$\epsilon>0$} 

\noindent
For $0< M< N$,
$$\log Z_{0,N}\geq \log Z_{N}(\{ \varphi_{M-1}=\varphi_{M} =0\})=\log Z_{0,M} +2\log \epsilon +\log Z_{M,N,\theta^{M}\omega}.$$
Therefore, $\{ \log Z_{0,N}+2\log \epsilon\}_{N\in \mathbb{N}}$ satisfies the ``super-additivity". The growth condition for $\mathbb{E}\log Z_{0,N}$
is given by the control of partition function. Let $l:= \#\{n: \varphi_n =0, 1\leq n\leq N-1\}$. Let $p\in \{0,1\}^{N}$. According to Lemma 3.1, the determinant for each path can be written as 
$$\sum_{|p|=N-1-l}\limits c_{ p} \;\prod_{i=0}^{N-1} \exp(\beta\omega_i)^{p_i}$$
where $\{c_p\}$ is a sequence of nonnegative integers. Notice that if $\beta=0$, the sum of $\{c_p\}$ is equal to the determinant in the nonrandom case. Let $L^{\varphi}$($L^{\varphi,\omega}$) be the matrix in the nonrandom(random) case. We have a equivalent expression and an upper bound for the partition function.
\[\begin{array}{rcl}
&&       \prod_{i=0}^{N}\exp(\frac{\beta}{2}\omega_i) Z^{\beta,\epsilon}_{N,\omega} \\
& =   &  \sum_{l=0}^{N-1} \epsilon^l (\sqrt{2\pi})^{-(l+1)} \prod_{i=0}^{N}\exp(\frac{\beta}{2}\omega_i)\left[\mbox{det}(L^{\varphi,\omega})\right]^{-1/2} \\
& =   &  \sum_{l=0}^{N-1} \epsilon^l (\sqrt{2\pi})^{-l}\left[\sum_{|p|=l+2}\limits c^{\varphi}_{ p} \;\prod_{i=0}^N \exp(-p_i \beta\omega_i)\right]^{-1/2}    \\
& =   &  \sum_{l=0}^{N-1} \epsilon^l (\sqrt{2\pi})^{-l} \left[\mbox{det}(L^{\varphi})\right]^{-1/2} \left[\sum_{|p|=l+2}\limits \frac{c^{\varphi}_{ p}}{\mbox{det}(L^{\varphi})} \;\prod_{i=0}^N \exp(-p_i \beta\omega_i)\right]^{-1/2}    \\

&\leq &  \sum_{l=0}^{N-1} \epsilon^l (\sqrt{2\pi})^{-l} \left[\mbox{det}(L^{\varphi})\right]^{-1/2} \left[\sum_{|p|=l+2}\limits \frac{c^{\varphi}_{ p}}{\mbox{det}(L^{\varphi})} \;\prod_{i=0}^N \exp(p_i \frac{\beta}{2}\omega_i)\right]   \\

\end{array}\]
Thus,
\[\begin{array}{rcl}
&&     \frac{1}{N}\mathbb{E}\log Z^{\beta,\epsilon}_{N,\omega} \\
&\leq &  \frac{1}{N}\mathbb{E}\log \left( \sum_{l=0}^{N-1} \epsilon^l (\sqrt{2\pi})^{-l} \left[\mbox{det}(L^{\varphi})\right]^{-1/2} \left[\sum_{|p|=l+2}\limits \frac{c^{\varphi}_{ p}}{\mbox{det}(L^{\varphi})} \;\prod_{i=0}^N \exp(p_i \frac{\beta}{2}\omega_i)\right] \right)  \\

&\leq &  \frac{1}{N}\log \left( \sum_{l=0}^{N-1} \epsilon^l (\sqrt{2\pi})^{-l} \left[\mbox{det}(L^{\varphi})\right]^{-1/2} \left[\sum_{|p|=l+2}\limits \frac{c^{\varphi}_{ p}}{\mbox{det}(L^{\varphi})} \;\mathbb{E}\prod_{i=0}^N \exp(p_i \frac{\beta}{2}\omega_i)\right] \right) \\

&=   &   \frac{1}{N}\log Z_{N}^{0,\epsilon M (\frac{\beta}{2})} +\frac{1}{N}\log (M(\frac{\beta}{2}))^{2}
\end{array}\]

We apply Liggett's version of subadditive ergodic theorem, $\frac{1}{N} \log Z_{0,N}$ converges a.s. and $\frac{1}{N}\mathbb{E}\log Z_{0,N}$ converges.

Again, we define the adjusted free energy for Laplacian model:
\begin{equation}
F_{N}(\beta,\epsilon):=\frac{1}{N}\mathbb{E}\log \mathcal{Z}^{\beta,\epsilon}_{N,\omega},\;\; F(\beta,\epsilon):=\lim_{N \rightarrow \infty}\limits F_{N}(\beta,\epsilon), \end{equation}
where
\begin{equation}
\mathcal{Z}^{\beta,\epsilon}_{N,\omega}:= \prod_{i=0}^{N} \exp(\frac{\beta}{2}\omega_i) Z^{\beta,\epsilon}_{N,\omega}\end{equation}

The quenched critical point is well-defined by 
\begin{equation}
\epsilon_c(\beta) :=\inf \{\epsilon : F(\beta,\epsilon)>0\}.
\end{equation}
\noindent
Similarly, we set the annealed critical point as
\begin{equation}\epsilon_c^a(\beta) :=\inf \{\epsilon : F^a(\beta,\epsilon)>0\}.\end{equation}

\subsection{Strong disorder regime: Iterated fractional moment estimates} 

\noindent
In this section, we know the exponent of the rate of the polynomial decay is $1$. This gives connections to the general pining model and the copolymer model which are discussed in \cite{Gia07,Gia11}. In this section, we prove the strong disorder regime based on the strategy mentioned in \cite{Gia11} Chapter 6, which is called the ``iterated fractional moment method''. The idea of this method is that for each $\beta$, finding a positive value $\Delta$ such that $F(\beta,\epsilon)=0$, where $\epsilon= \epsilon_c^a(\beta)e^{\Delta}$. One observation is that
\begin{equation}
F(\beta,\epsilon) =\lim_{N\rightarrow \infty}\limits \frac{1}{N} \mathbb{E}\log \mathcal{Z}_N
=\lim_{N\rightarrow \infty}\limits \frac{1}{\gamma N} \mathbb{E}\log \mathcal{Z}_N^{\gamma}\leq \lim_{N\rightarrow \infty}\limits \frac{1}{\gamma N} \log \mathbb{E} \mathcal{Z}_N^{\gamma}
\end{equation}
for any $\gamma>0$. Since the annealed quantity $\mathbb{E} \mathcal{Z}_N^{\gamma}$ is more tractable, we will choose $\gamma$ and $\Delta$, such that $\mathbb{E} \mathcal{Z}_N^{\gamma}$ is bounded by a constant. Thus, $F(\beta,\epsilon)=0$, and $\log \epsilon_c(\beta)-\log \epsilon_c^a(\beta)\geq \Delta >0$.\\

The following classical result helps the fractional moment estimate.\\

\begin{lemma} (\cite{HLP67} Chapter 2.1) Let $0< \gamma <1$ and $\{a_n\}_n$ is a positive sequence. Then
$$(a_1+\cdots+a_n)^{\gamma} < a_1^{\gamma} +\cdots+ a_n^{\gamma}.$$
\end{lemma}

If for given $\beta$ and $\Delta$ we can find a fixed number $k$ and $\gamma\in(0,1)$ such that
\begin{equation}
\rho:= (e^{\Delta}R(\beta)^{-1})^{\gamma} \sum^{\infty}_{n=k+1} \sum_{s=0}^k K(n-s)^{\gamma}  \mathbb{E}e^{\gamma \psi(\beta\bar{\omega}[0,n-s))} A_{s} \leq 1,
\end{equation}
then we have
\begin{equation}
A_N\leq\rho \max\{A_0,...,A_{N-k-1}\}
\end{equation}
for $N>k$, which implies that $A_N \leq \max\{A_0,...,A_{k}\}$ and hence $F(\beta,\epsilon_c^a(\beta)e^{\Delta})=0$, that is, $\log \epsilon_c(\beta)-\log \epsilon_c^a(\beta)\geq \Delta$. \\

The proof of the following proposition is based on \cite{Gia11}.\\

For $k-R< s\leq k$, 
we introduce the "tilting measure" $\tilde{\mathbb{P}} :=  \tilde{\mathbb{P}}_{n,\lambda}$ for $n\in \mathbb{N}$, $\lambda \in\mathbb{R}$ and
\begin{equation}
\frac{\mbox{d}  \tilde{\mathbb{P}}_{n,\lambda}}{ \mbox{d} \mathbb{P}} (\omega):= \frac{1}{M(-\lambda)^n} \exp(- \lambda\sum_{i=0}^{n-1} \omega_i).\end{equation}
Now, we use the H\"{o}lder's inequality
$$ A_s = \tilde{\mathbb{E}} \left[ (Z_s)^{\gamma}  \frac{d\mathbb{P}}{d\tilde{\mathbb{P}} }\right] \leq  \left( \tilde{\mathbb{E}}  \left[ \frac{d\mathbb{P}}{d\tilde{\mathbb{P}} }\right]^{1/(1-\gamma)}  \right)^{1-\gamma} \left( \tilde{\mathbb{E}} Z_s \right)^{\gamma} $$
For the first term, we have 
$$\left( \tilde{\mathbb{E}}  \left[ \frac{d\mathbb{P}}{d\tilde{\mathbb{P}} }\right]^{1/(1-\gamma)}  \right)^{1-\gamma}
=\exp\left(  \gamma s \log M(-\lambda)+ (1-\gamma)s \log M\left(\lambda \frac{\gamma}{1-\gamma }\right) \right)$$
$$\leq \exp \left( \gamma C_M(c\beta^2 k)+ (1-\gamma) C_M (c\beta^2 k)\left(\frac{\gamma}{1-\gamma} \right)^2\right)  \leq \exp\left( \frac{C_M}{C(\beta)}\frac{ \gamma}{1-\gamma}\right)  $$
where $2C_M:=  \mbox{max}_{|t|\leq 1}(\log M(t))''$ and provided the arguments of $M$ are less than 1 by choosing $c$ small.

\noindent
From $\check{Z}_{0,n}^{\epsilon}   \sim \frac{C_{\epsilon}}{\epsilon^2 n^2}$, we know this is close to the case $\alpha=1$ in \cite{Gia11}. The proof is delicate in this case. Here, we sketch the proof, one can see details in \cite{Gia11} Chapter 6. Denote $A_N:=\mathbb{E}\mathcal{Z}_{N,\omega}^{\gamma}$. By using Lemma 2.5, we have for $N>k$
\[\begin{array}{rcl}
A_N  &\leq  &  \epsilon^{2\gamma} \sum_{n=k+1}^N A_{N-n} \sum_{s=0}^{k}\mathbb{E}\check{\mathcal{Z}}^{\gamma}_{n-s}  \;A_{s}\\

     &\leq  &  \epsilon^{2\gamma} \sum_{n=k+1}^N A_{N-n} \sum_{s=0}^{k} (\mathbb{E}\check{\mathcal{Z}}_{n-s})^{\gamma}  \;A_{s}\\
     
     &\leq  &  \epsilon^{2\gamma} \sum_{n=k+1}^N A_{N-n} \sum_{s=0}^{k} (\frac{C_{\beta}}{(n-s+1)^2})^{\gamma}  \;A_{s}\\
\end{array}   \]

For given $\beta$ and $\Delta$ we try to find $k$ and $\gamma\in(0,1)$ such that
$$\rho:= \epsilon^{2\gamma} \sum_{n=k+1}^{\infty} \sum_{s=0}^{k} \left(\frac{C_{\beta}}{(n-s+1)^2}\right)^{\gamma}  \;A_{s} $$
is small.\\

\noindent
{\bf Proof of Proposition 1.5.}\\
\noindent
First, $\gamma= \gamma(k)=1-(1/\log k)$. As suggested in \cite{Gia11} we choose 
$$  \Delta:=\frac{ c\beta^2}{( \log(1+\frac{1}{\beta}) )^2} \; , \; k:= \lfloor\frac{ ( \log (1+\frac{1}{\beta}) )^2 }{c\beta^2} \rfloor 
\;\mbox{and} \;   \lambda :=\frac{\sqrt{c}\beta}{( \log(1+\frac{1}{\beta}) )^2}$$ 
Notice that
$$A_s\leq (\mathbb{E}\mathcal{Z}_s)^{\gamma}=\left[\frac{\exp(sF^{V_{\beta}}(\Delta)) P_{\Delta}(s+1\in \chi) }{\epsilon_c^a(\beta)^2\exp(2\Delta)}\right]^{\gamma}
\leq \left[\frac{\exp(kF^{V_{\beta}}(c\beta^2/\log^2(1+1/\beta))) }{\epsilon_c^a(\beta)^2}\right]^{\gamma} $$
is bounded for $s\leq k$.   
We estimate
\[\begin{array}{rcl}
\tilde{E}\mathcal{Z}_s  &=    & \mathbb{E}(\mathcal{Z}_s\frac{\exp(-\lambda \sum_{i=0}^s \omega_i )}{M(-\lambda)^{s+1}})\\
              &=    & Z(s,\epsilon_c^a(\beta)\exp(-\beta\lambda/2)\exp(\Delta),V_{\beta})\cdot \exp(-\beta\lambda)\\ 
              &=    & E_{s,\epsilon_c^a(\beta)}\left( e^{L_s(\Delta-\beta\lambda/2)}\right) \cdot Z(s,\epsilon_c^a(\beta),V_{\beta}) \cdot \exp(-\beta\lambda)\\
              &=    & \mathcal{E}_{\epsilon_c^a(\beta)}\left( e^{L_s(\Delta-\beta\lambda/2)}  |s+1\in\chi \right) \cdot \frac{\mathcal{P}_{\epsilon_c^a(\beta)}(s+1\in \chi)}{\epsilon_c^a(\beta)^2} \cdot \exp(-\beta\lambda)\\
              &=    & \mathcal{E}_{\epsilon_c^a(\beta)}\left( e^{L_s(\Delta-\beta\lambda/2)}  1_{\{s+1\in\chi\}} \right) \cdot \frac{1}{\epsilon_c^a(\beta)^2} \cdot \exp(-\beta\lambda)
\end{array}\]
The second equality is due to the property of Gaussian variables. Recall that the quantity $L_s$ is the cardinality of $\{ 0 < n \leq s : \varphi_n =0\}$, and $i_s$ is the cardinality of $\{ 0 < n \leq s : \varphi_{n-1} =\varphi_n =0\}$, thus, $L_s\geq i_s$. By Proposition 3.3, the double-return sequence $\{\chi_k\}_{k\geq 0}$ is a genuine renewal process under $\mathcal{P}_{\epsilon_c^a(\beta)}$ with renewal distribution 
$$K(n) \sim \frac{C(\epsilon_c^a(\beta))}{n^2}.$$ 
Based on the value of $\Delta$ and $\lambda$ we choose, $\Delta -\beta\lambda/2<0$ and
$$\tilde{E}\mathcal{Z}_s \leq  \mathcal{E}_{\epsilon_c^a(\beta)}\left( e^{i_s(\Delta-\beta\lambda/2)} \right) \cdot \frac{1}{\epsilon_c^a(\beta)^2}. $$

The rest of proof goes the same as \cite{Gia11} Chapter 6, we get $\tilde{E}\mathcal{Z}_s$ arbitrarily small if $c$ is small.

Remark. For general charges, the estimate of the tilted partition function would be
\[\begin{array}{rcl}
\tilde{E}\mathcal{Z}_s  &=    & \mathbb{E}(\mathcal{Z}_s\frac{\exp(-\lambda \sum_{i=0}^s \omega_i )}{M(-\lambda)^{s+1}})\\
              &=    & Z(s,\epsilon_c^a(\beta)\exp(\Delta),V_{\beta,-\lambda})\\ 
              &=    & E_{s,\epsilon_c^a(\beta,-\lambda)}\left( \left(  \frac{\epsilon_c^a(\beta)\exp(\Delta)}{\epsilon_c^a(\beta,-\lambda)} \right)^{ L_s}\right) \cdot Z(s,\epsilon_c^a(\beta,-\lambda),V_{\beta,-\lambda})\\
              &=    & \mathcal{E}_{\epsilon_c^a(\beta,-\lambda)}\left( \left(  \frac{\epsilon_c^a(\beta)\exp(\Delta)}{\epsilon_c^a(\beta,-\lambda)} \right)^{L_s}|s+1\in\chi \right)\cdot \frac{\mathcal{P}_{\epsilon_c^a(\beta,-\lambda)}(s+1\in \chi)}{\epsilon_c^a(\beta,-\lambda)^2}\\
              &\leq & \mathcal{E}_{\epsilon_c^a(\beta,-\lambda)}\left( \left(  \frac{\epsilon_c^a(\beta)\exp(\Delta)}{\epsilon_c^a(\beta,-\lambda)} \right)^{L_s} \right)\cdot \frac{1}{\epsilon_c^a(\beta,-\lambda)^2}
\end{array}\]

Based on the fact $\epsilon_c^a(\beta) \leq \epsilon_c^a(\beta,-\lambda)$, we get
$$\tilde{E}\mathcal{Z}_s \leq \mathcal{E}_{\epsilon_c^a(\beta,-\lambda)}\left( \left(  \frac{\epsilon_c^a(\beta)\exp(\Delta)}{\epsilon_c^a(\beta,-\lambda)} \right)^{L_s} \right) \cdot \frac{1}{\epsilon_c^a(\beta)^2} $$

However, it's not obvious that there exists a positive constant C, such that
$$ \log \frac{\epsilon_c^a(\beta)}{\epsilon_c^a(\beta,-\lambda)} \leq - C \beta\lambda .$$

\section{Special determinants} 

\subsection{Proof of Lemma 2.4} 

\noindent
Given a positive sequence $(b_0,...,b_n)$, $B_{n-1\times n-1}$ is a symmetric matrix and its upper triangle part is defined as following:
\[B_{ij}= \left \{ \begin{array}{ll}
b_{i-1}+4b_i+b_{i+1}, & i=j\\
-2b_i-2b_{i+1},       & i=j-1\\
b_{i+1},              & i=j-2 \\
0,                    & i\leq j-3.
\end{array}\right.\]
For example,
\[ B_{5\times 5} = \left[ \begin{array}{ccccc}  
b_0+4b_1+b_2 & -2b_1-2b_2   & b_2          & 0						& 0 \\
-2b_1-2b_2   & b_1+4b_2+b_3 & -2b_2-2b_3   & b_3					& 0 \\
b_2          & -2b_2-2b_3   & b_2+4b_3+b_4 & -2b_3-2b_4	  & b_4 \\
0            &b_3           & -2b_3-2b_4   & b_3+4b_4+b_5 & -2b_4-2b_5 \\
0            &0             &b_4           & -2b_4-2b_5   & b_4+4b_5+b_6
\end{array} \right]
\]
$B_{n-1\times n-1}$ is positive-definite, for $\varphi^tB\varphi= \sum_{m=0}^{n} b_m(\Delta \varphi_m)^2\geq 0$. Let $D(n-1):= \sum_{k=1}^n \sum_{i=0}^{n-k} k^2b_i^{-1}b_{i+k}^{-1}$.
We claim that the determinant of $B_{n-1}$ is $\prod_{i=0}^n b_i\cdot D(n-1)$. The proof is given by row operations and the mathematical induction. We use $B_5$ to elaborate the ideas. First, let new rows be $r'_i=\sum_{j=1}^i (i-j+1)r_j\;\; i=1,...,5$. For the new matrix, add twice of the second column to the first one. Then we have a matrix having the same determinant as $B_5$:
\[ \left[ \begin{array}{ccccc}  
b_0-3b_2    & -2b_1-2b_2   & b_2      & 0						& 0 \\
2b_0+2b_3   & -3b_1+b_3    & -2b_3    & b_3					& 0 \\
3b_0        & -4b_1        & b_4      & -2b_4	  & b_4 \\
4b_0        & -5b_1        & 0        & b_5     & -2b_5 \\
5b_0        & -6b_1        & 0        & 0       & b_6
\end{array} \right]
\]
Grab the common factor $b_0$ and $b_{1}$ from colume 1 and colume 2, respectively. Also, grab the common factor $b_{i+1}$ from the ith row. It suffices to show that the deteminant of 
\[ B'=\left[ \begin{array}{ccccc}  
b_2^{-1}-3b_0^{-1}    & -2b_2^{-1}-2b_1^{-1}   & 1   & 0	& 0 \\
2b_3^{-1}+2b_0^{-1}   & -3b_3^{-1}+b_1^{-1}    & -2  & 1	& 0 \\
3b_4^{-1}             & -4b_4^{-1}             & 1   & -2	& 1 \\
4b_5^{-1}             & -5b_5^{-1}             & 0   & 1  & -2 \\
5b_6^{-1}             & -6b_6^{-1}             & 0   & 0  & 1
\end{array} \right]
\]
is equal to D(5), which is $D(4)+ b_6^{-1}(\sum_{j=0}^5 (6-j)^2 b_j^{-1})$. Now, we expand the determinant by the last row, and notice that the determinant of the principle $4\times 4$ matrix is D(4). So it remains to show that
\[ (-1)^4\cdot 5\left| \begin{array}{cccc}  
-2b_2^{-1}-2b_1^{-1}   & 1   & 0	& 0 \\
-3b_3^{-1}+b_1^{-1}    & -2  & 1	& 0 \\
-4b_4^{-1}             & 1   & -2	& 1 \\
-5b_5^{-1}             & 0   & 1  & -2 \\
\end{array} \right|+
(-1)^5\cdot (-6)\left| \begin{array}{cccc}  
b_2^{-1}-3b_0^{-1}    & 1   & 0	& 0 \\
2b_3^{-1}+2b_0^{-1}   & -2  & 1	& 0 \\
3b_4^{-1}             & 1   & -2	& 1 \\
4b_5^{-1}             & 0   & 1  & -2 \\
\end{array} \right|
\]
is equal to $\sum_{j=0}^5 (6-j)^2 b_j^{-1}$. For $B_{n-1}$, after we follow the same procedure, it suffices to show that $\sum_{j=0}^{n-1} (n-j)^2 b_j^{-1}$ equals the determinant of a $(n-2)\times(n-2)$ matrix
\[ (-1)^{n-2} \left| \begin{array}{ccccc}  
-3nb_0^{-1}-2(n-1)b_1^{-1}-(n-2)b_2^{-1}    & 1   & 0   & 0\cdots 0	& 0 \\
2nb_0^{-1}+(n-1)b_1^{-1}-(n-3)b_3^{-1}      & -2  & 1	  & 0\cdots 0  & 0   \\
-(n-4)b_4^{-1}                              & 1   & -2  & 1\cdots 0	& 0  \\
\vdots                                      &     &     &           &        \\
-b_{n-1}^{-1}                               & 0   &0    & 0\cdots 1  &-2  
\end{array} \right|
\]
which is the same as
\[ \left| \begin{array}{ccccc}  
3nb_0^{-1}+2(n-1)b_1^{-1}+(n-2)b_2^{-1}    & -1   & 0   & 0\cdots 0  	& 0 \\
-2nb_0^{-1}-(n-1)b_1^{-1}+(n-3)b_3^{-1}    & 2    & -1	& 0\cdots 0    & 0   \\
(n-4)b_4^{-1}                              & -1   & 2   & -1\cdots 0   & 0  \\
\vdots                                     &      &     &             &        \\
b_{n-1}^{-1}                               & 0    &0    &  0\cdots -1 &2  
\end{array} \right|
\]
Notice that the right bottom is the matrix $A_{(n-3)\times(n-3)}$ with $a_i=1\;\forall i$. The proof is completed by expanding the determinant by the first column. Again, notice that every term in the det($B_{n-1}$) is of degree $(n-1)$ and has no multiplicity.

\subsection{A more general case to Lemma 2.4.} 

\noindent
For general cases, if the path $\{\varphi_n\}_{n\leq N+1}$ hits 0 between 0 and N, we still can compute the corresponding determinant by deleting the $m$th column and $m$th row if $\varphi_m=0$. For example, if $N=6$ and only $\varphi_4=0$, the underlying matrix is
\[\left[ \begin{array}{cccc}  
b_0+4b_1+b_2 & -2b_1-2b_2   & b_2          	& 0 \\
-2b_1-2b_2   & b_1+4b_2+b_3 & -2b_2-2b_3   	& 0 \\
b_2          & -2b_2-2b_3   & b_2+4b_3+b_4  & b_4 \\
0            &0             &b_4            & b_4+4b_5+b_6
\end{array} \right]\]

It's natural to guess that every term in the determinant is of degree 4 and has no multiplicity.\\

\begin{lemma}
Given a path $\{\varphi_n\}_{n\leq N+1}$ and $r=\#\{n: \varphi_n =0, 1\leq n\leq N-1\}\geq 1$. Every term in the corresponding determinant is of degree $(N-1-r)$ and has no multiplicity.
\end{lemma}

\noindent
Proof. We prove it by induction. Given a path $\{\varphi_n\}_{n\leq N+2}$, we need to show that the degree is $(N-r)$. Let $m=\sup \{n: \varphi_n =0, 1\leq n\leq N\}$.  If $m=1$, by the previous lemma, every term in the determinant is of degree $(N-1)$. Note that $r=1$ since $m=1$. If $m=N$, by the induction hypothesis, every term in the determinant is of degree $(N-1-(r-1))$. If $2\leq m\leq N-1$ and $\varphi_{m-1}=0$, the corresponding matrix becomes
\[\left[ \begin{array}{cc}  
A & 0  \\
0 & C 
\end{array} \right]\]
where A is a $[(m-2)-(r-2)]\times [(m-2)-(r-2)]$ matrix, and C is a $(N-m)\times (N-m)$ matrix. Thus, the determinant is $\mbox{det}(A)\mbox{det}(C)$. By the induction hypothesis and previous lemma, every term in $\mbox{det}(A)$ is of degree $m-r$, and every term in $\mbox{det}(c)$ is of degree $N-m$. On the other hand, if $\varphi_{m-1}\neq 0$, the corresponding matrix is still positive-definite and can be written as 
\[\left[ \begin{array}{cc}  
A & E  \\
E^*   & C 
\end{array} \right]\]
where A is a $[(m-1)-(r-1)]\times [(m-1)-(r-1)]$ matrix, and C is a $(N-m)\times (N-m)$ matrix. Clearly, the only nonzero term in $E$ is $e_{m-r,m-r+1}=b_m$. Set $X=-A^{-1}E$, the determinant is equal to
\[\mbox{det}\left[ \begin{array}{cc}  
I & 0  \\
X^* & I 
\end{array} \right]
\left[ \begin{array}{cc}  
A & E  \\
E^* & C 
\end{array} \right]
\left[ \begin{array}{cc}  
I & X  \\
0 & I 
\end{array} \right]
=\mbox{det}\left[ \begin{array}{cc}  
A & 0  \\
0 & C-E^*A^{-1}E 
\end{array} \right]\]
Let $A_{-1}$ be the matrix deleting the last column and row from A, and $C_{-1}$ be the matrix deleting the first column and row from C. (If A is of dimension 1, let $A_{-1}=I$, so is $C_{-1}$). Let $A'=A|_{b_m=0}$ and $C'=C|_{b_m=0}$. We have 
$$\mbox{det}(A) =  \mbox{det}(A')+b_m\mbox{det}(A_{-1}),$$
$$\mbox{det}(C-E^*A^{-1}E) = \mbox{det}(C')+[b_m-\frac{\mbox{det}(A_{-1})}{\mbox{det}(A)}b_m^2]\mbox{det}(C_{-1}).$$
So the underlying determinant is equal to
\[\begin{array}{rcl}  
\mbox{det}(A)\mbox{det}(C-E^*A^{-1}E) & =&  \mbox{det}(A)\mbox{det}(C')  \\
&    &   +[\mbox{det}(A')+b_m\mbox{det}(A_{-1})][b_m-\frac{\mbox{det}(A_{-1})}{\mbox{det}(A)}b_m^2]\mbox{det}(C_{-1})\\
&=   & \mbox{det}(A)\mbox{det}(C') \\
&    &   +[\mbox{det}(A')b_m+b_m^2\mbox{det}(A_{-1})-\mbox{det}(A_{-1})b_m^2]\mbox{det}(C_{-1})\\
&=   & \mbox{det}(A)\mbox{det}(C')+ \mbox{det}(A')b_m\mbox{det}(C_{-1})
\end{array} \]
The degree in the first term is $(m-r)+(N-m)$, and the degree in the second term is $(m-r)+1+(N-m-1)$. It's easy to see that there is no multiplicity, which ends the proof.

\end{document}